\documentclass[twoside,10pt,reqno]{amsart}
\usepackage{amsmath}
\usepackage{amssymb}
\usepackage[initials,alphabetic]{amsrefs}
\newtheorem{thm}{Theorem}[section]
\newtheorem{prop}[thm]{Proposition}
\newtheorem{defn}[thm]{Definition}
\newtheorem{lemma}[thm]{Lemma}
\newtheorem{problem}[thm]{Problem}

\newtheorem{remark}[thm]{Remark}

\begin{document}

\title[Yamabe problem on noncompact manifolds with boundary]{
The Zero Scalar Curvature Yamabe problem on noncompact 
manifolds with boundary}
\author{Fernando Schwartz}
\email{fernando@math.duke.edu}
\address{Cornell University, Ithaca, NY, and Stanford University, Stanford, CA}
\curraddr{Duke University, Durham, NC}
%\date{June 2004}
\begin{abstract}
Let $(M^n,g),~n\ge 3$ be a noncompact complete Riemannian manifold with
compact boundary
and $f$  a smooth function on $\partial M$. 
In this paper we show that for a large class of such manifolds, there exists a 
metric within the conformal class of $g$ that is complete, has zero scalar 
curvature on $M$ and has mean curvature $f$ on the boundary.\\
The problem is equivalent to finding a positive solution
to an elliptic equation with a non-linear boundary condition with critical 
Sobolev exponent.
\end{abstract}

\dedicatory{Dedicated to the memory of Prof. Jos\'e F. Escobar}
\thanks{The author would like to thank Prof. Richard Schoen for his thorough support
and great comments. This work was supported in part by NSF grant \# 
DMS-0223098}

\maketitle

\section{Introduction}

The celebrated Riemann Mapping Theorem states that any simply connected 
region in the plane is conformally diffeomorphic to a disk. This theorem
is less successful in higher dimensions since very few domains are conformally 
diffeomorphic to the ball. Nevertheless, we can still ask whether a
manifold with boundary is conformally diffeomorphic to a manifold that resembles the
ball, namely to one that has zero scalar curvature and constant mean curvature on 
its boundary. Escobar studied this problem in \cite{E92}. He showed that
most compact manifolds with boundary admit such conformally related metrics.

A generalization of this problem is the so-called {\it prescribed mean curvature
problem}. Let $(M^n,g),n\ge 3$ be a manifold with boundary and $f\in C^\infty(\partial M)$.

\begin{problem}\label{yamb}
Does there exist a metric conformally equivalent
to $g$ that is complete, scalar flat and has mean curvature $f$ on $\partial M$?
\end{problem}
Escobar and Garcia \cite{EG04} studied this problem 
on $(B^3,\delta_{ij})$. They  proved that a Morse function 
is the mean curvature of a scalar-flat metric $g\in [\delta_{ij}]$ 
if it satisfies some Morse inequalities. They paralleled Schoen and 
Zhang's \cite{SZ96} blow-up analysis for the prescribed scalar curvature problem
on $S^3$. In both cases, a general solution is unexpected because of 
the Kazdan-Warner obstruction \cite{KW75a}. (See
 \cite{KW74} for the prescribed Gaussian curvature problem
on open 2-manifolds.)

In this paper we address Problem \ref{yamb} on a large class of {\bf noncompact}
manifolds with boundary $(M^n,g),~n\ge3$. As a corollary of
Theorem \ref{thm1} about PDEs we get:
\begin{thm} \label{bb}
   Any smooth function $f$ on $\partial M$ can be realized as the mean curvature 
   of a complete scalar flat metric conformal to $g$.
\end{thm}
In contrast with the compact case, no topological obstructions on $f$
arise. This is a surprising phenomena.

\section{Preliminaires}

Let $(M^n,g),~n\ge 3$ be a complete, $n$-dimensional Riemannian manifold with
boundary $\partial M\neq \emptyset$.
Denote by $\tilde{g}=u^{4/(n-2)}g$ a metric conformally related to $g$, where $u>0$
is a smooth function. 

It is a standard fact that the relation between the scalar curvature $R(g)$ of the metric $g$ 
and the scalar curvature $R(\tilde{g})$ of the metric $\tilde{g}$ is given by

\begin{eqnarray}\label{eq1}
R(\tilde{g})=-\frac{4(n-1)}{n-2}\frac{L_gu}{u^{(n+2)/(n-2)}},
\end{eqnarray}
where $L_g=\Delta_g-\frac{n-2}{4(n-1)}R(g)$, and $\Delta_g$ is the Laplacian 
calculated with respect to the metric $g$.

The relation between the mean curvature of the boundary $h(g)$ of the metric 
$g$, and the mean curvature of the boundary $h(\tilde{g})$ of the
metric $\tilde{g}$ is given by

\begin{eqnarray}\label{eq2}
h(\tilde{g})=\frac{2}{n-2}\frac{B_gu}{u^{n/(n-2)}},
\end{eqnarray}
where $B_g=\frac{\partial}{\partial\eta}+\frac{n-2}{2}h(g)$ and $\partial /\partial\eta$ 
is the outward-pointing normal derivative on $\partial M$ calculated with respect to 
the metric $g$.

\begin{remark}
The exponent $n/(n-2)$ of  equation (\ref{eq2})  is called a {\bf critical exponent}
since the Sobolev trace embedding $H^1(M)\hookrightarrow L^{q}(\partial M)$ ceases 
to be compact for $q\ge 2(n-1)/(n-2)$. This condition rules out the direct method of
minimization to prove existence of solutions.
\end{remark}

It follows directly from the above discussion that
finding a conformally related metric $\tilde{g}=u^{4/(n-2)}g$ on $M$
that is scalar flat (i.e. has zero scalar curvature) and has prescribed 
mean curvature $f$ on the boundary is equivalent to finding smooth 
$u>0$ on $M$ that satisfies equation (\ref{eq1}) with $R(\tilde{g})\equiv 0$
and equation (\ref{eq2}) with $h(\tilde{g})\equiv f$.

In this paper we find such $u$ for a more general problem, the so-called
{\bf supercritical} equation, in which the critical exponent $n/(n-2)$ of  
(\ref{eq2}) is replaced by an arbitrary number $\beta>1$.

\begin{defn}
Let $(M^n,g)$ be a complete, noncompact Riemannian manifold.
On each end $E$ of $M$, consider the volume of the set obtained by intersecting
$E$ with the geodesic ball of radius $t$ centered at some fixed $p\in M$, and
denote it by $V_E(t)$. We say that the end $E$ is {\bf large} if
\[ \int_1^\infty \frac{t}{V_E(t)} dt <\infty.\]
\end{defn}

Suppose that the Ricci curvature of $M$ satisfies
$Ric_M(x)\ge -(n-1)K(1+r(x))^{-2}$, where  $K\ge0$ is some constant and 
$r(x)$ is the distance from $x$ to some fixed point $p$. By Li and Tam's \cite{LT95} paper,
on any large end $E$ of $M$ there exists 
a harmonic, non-negative function $v_E$ (a \it barrier\rm), which is asymptotic to 1 on $E$
and it is exactly zero on the boundary of a large ball intersected with the end.

Throughout this paper $M$ will be a manifold that satisfies the above  bound on the Ricci tensor.

\begin{defn}\label{positive}
We say that $(M,g)$ is {\bf positive} if it is complete, scalar flat, and has positive 
mean curvature on the boundary.
\end{defn}

\begin{remark}
If $(M,g)$ is positive it has a positive first 
eigenvalue for the following problem:
\begin{eqnarray*}\label{CN} 
 \left\{
\begin{array}{lllr}
L_{g} u & = & 0  &\mbox{in}~M,\\
B_{g} u & = & \lambda u   & \mbox{on}~\partial M.
\end{array}
\right.
\end{eqnarray*}
Conversely, if the first eigenvalue of the above problem
is positive, then $(M,g)$ admits a conformally-related 
sclar flat metric $g'$ that has postive mean curvature on the
boundary, but this metric {\bf  may not be complete}. Theorem \ref{thm1}
still applies for $(M,g')$ provided it has large ends, since completeness
is not used in the proof.
Nevertheless, the new metric $\tilde{g}=u^{4/(n-2)}g'$, which is scalar flat and
has prescribed mean curvature $f$ on the boundary, may also be incomplete. 
\end{remark}

\begin{thm}\label{thm1}
Let $(M,g)$ be a noncompact positive Riemannian manifold with compact
boundary and finitely many ends, all of them large. Let $f$ be a smooth function on 
$\partial M$ and $\beta>1.$
There exists $\epsilon,\delta>0$ and a smooth function 
$\epsilon\le u\le \epsilon+\delta$ on $M$ with
\begin{eqnarray}\label{cf}
\left\{
\begin{array}{lllr}
L_gu  & = & 0 &\mbox{in}~M,\\
B_g u &  = & \frac{n-2}{2}f u^{\beta} & \mbox{on}~\partial M.
\end{array}
\right.
\end{eqnarray}
When $\beta$ is the critical exponent $n/(n-2),~\tilde{g}=u^{4/(n-2)}g$
is a complete, scalar flat metric on $M$, with mean
curvature $h(\tilde{g})\equiv f$.
\end{thm}

\begin{remark}
For the $\beta=n/(n-2)$ case, the bound  $\epsilon\le u\le \epsilon+\delta$ 
guarantees a complete metric $\tilde{g}=u^{4/(n-2)}g$.
\end{remark}

\begin{remark}
Since $(M,g)$ is positive we have that $L_g\equiv \Delta_g.$
\end{remark}

A very important class of examples of positive {\it noncompact} manifolds 
with boundary is obtained by removing submanifolds of large codimension out 
of positive {\it compact} manifolds with boundary. We refer the reader to the Appendix 
for more details on the construction.

The proof of Theorem \ref{thm1} is divided in two parts. In
Section \ref{M} we prove that an iterative process using
sub- and super-solutions converges to a solution 
of (\ref{cf}).
In Section \ref{E} we construct the appropriate sub- and super-solutions.
Theorem \ref{bb} follows by choosing $\beta=n/(n-2)$ in Theorem \ref{thm1}.

%------------------------------------------------------------------------------------------------------------

\section{Method of sub- and super-solutions}
\label{M}
%------------------------------------------------------------------------------------------------------------

In this section we adapt a method of sub- and super-solutions
to our setting. (See \cite{KW75b} for general properties of  sub- 
and super-solution methods on semilinear elliptic problems.)
We begin by proving a form of maximum principle
on a compact piece of $M$ that contains $\partial M$.

Let $u\in C^{2}(M)\cap C^1(\bar{M})$, and define the operators
\[
\begin{array}{lllr}
L_\lambda u & := & \Delta_g u -\lambda u & \mbox{in}~M,\\
B_\gamma u & := &\frac{\partial u}{\partial \eta}+(\frac{n-2}{2}h_g+\gamma) u  & 
\mbox{on}~\partial M,
\end{array}
\]
for $\lambda,\gamma>0$ 
fixed large numbers.

\begin{prop}[Maximum Principle] \label{mp} Let $M_1\subseteq M$ be
a compact piece of $M$ containing $\partial M$, with smooth boundary $\partial M_1 = \partial M \cup N$.
Suppose that  $u\in C^{2}(M_1)\cap C^1(\bar{M_1})$ satisfies:
\[ \left\{
\begin{array}{lllr} 
L_\lambda u& \geq &  0& \mbox{in}~ M_1,\\
B_\gamma u &\leq & 0&\mbox{on} ~\partial M,\\
u & \le& 0 &\mbox{on}~ N.
\end{array}
\right.\]
 Then $u\leq 0$ on $ M_1$.
\end{prop}

\begin{proof} Put $w(x)=\max \{0,u(x)\}$, so that $w=0$ on $N$. Recall that $\min_{\partial M}
 h_g >0$. We get:
\begin{eqnarray*}
0 & \leq & \int_{M_1} (L_\lambda u)w - \int_{\partial M}(B_\gamma u)w\\
&= & -\int_{M_1}\nabla u\cdot \nabla w - \lambda \int_{{M_1}} uw - \gamma
\int_{\partial M} uw\\
& = & -\int_{M_1} |\nabla w|^2-\lambda \int_{{M_1}} w^2-\gamma\int_{\partial M}w^2.
\end{eqnarray*}
Hence $w=0$ in $M_1$, and so $u\leq 0$ in $M_1$.
\end{proof}

\begin{defn}
A   sub-solution (resp. super-solution) of equation (\ref{cf}) is a
function $u_- $ (resp. $u^+$) in $C^{2}(M)\cap C^1(\bar{M})$ with
\[
\left\{
\begin{array}{lllr}
\Delta_g u_- & \geq & 0 & \mbox{in}~M,\\
B_g u_-- \frac{n-2}{2} f u_-^\beta & \leq  & 0 & 
\mbox{on}~\partial M,
\end{array}
\right.
\]
respectively
\[
\left\{
\begin{array}{lllr}
\Delta_g u^+ & \leq & 0 & \mbox{in}~M,\\
B_g u^+- \frac{n-2}{2} f (u^+)^\beta & \geq  & 0 & 
\mbox{on}~\partial M.
\end{array}
\right.
\]

\end{defn}

\begin{thm} 
\label{main}
If there exist sub- and super-solutions $u_-,u^+\in C^\infty(M)$ with
$0\le u_-\le~u^+\le c_0$, then there exists
a smooth solution $u$ of equation (\ref{cf}) with $u_-\leq u \leq u^+$.
\end{thm}

\begin{proof} We will show that the statement holds for all
compact pieces $M_1\subseteq M$ as above. Then we take
pieces converging to $M$ and construct a global solution.

Let $M_1$ be a compact neighborhood of $\partial M$ in $M$
with smooth boundary $\partial M_1=\partial M \cup N$. Let 
$\lambda,\gamma>0$ be large enough so that (\ref{induct})
admits a solution. Let $u_0=u^+|_{M_1}$, 
and define inductively $u_i\in C^{2}(M_1)\cap C^1(\bar{M_1})$,
$i=1,2,\dots$, to be the unique solution to
\begin{eqnarray} \label{induct}
\left\{
\begin{array}{lllr}
L_\lambda u_i &  = & -\lambda u_{i-1} & \mbox{in} ~M_1,\\
B_\gamma u_i &  = &\frac{n-2}{2} f u_{i-1}^\beta + \gamma u_{i-1} & \mbox{on}~\partial M,\\
u_i & = & u_{i-1} & \mbox{on}~N.
\end{array}
\right.
\end{eqnarray}

\it Claim.
\rm
We have $u_-\leq \cdots \leq u_i\leq u_{i-1}\leq \cdots \leq u^+$.

To prove the claim, we will use induction twice. First, to show that the
sequence $\{u_i\}$ is non-increasing and bounded by $u^+$. Then, 
to prove that it is bounded below by $u_-$.
 
To check the first induction step, we see that $L_\lambda (u_1-u_0) = 
(\Delta u_1-\lambda u_1)-(\Delta u_0
-\lambda u_0) = -\lambda u_0 - \Delta u_0 +\lambda u_0 = -\Delta u_0
\geq 0$, because $u_0=u^+$ is a super-solution.

On the other hand, one has
 
\begin{eqnarray*}
B_\gamma
 (u_1-u_0)& =& \frac{\partial u_1}{\partial \eta} + (\frac{n-2}{2}h_g +\gamma) u_1
  -\frac{\partial u_0}{\partial
 \eta} - (\frac{n-2}{2}h_g +\gamma)u_0\\
& =&  \frac{n-2}{2}fu_0^\beta +\gamma u_0 -\frac{\partial u_0}{\partial
 \eta}  - (\frac{n-2}{2}h_g +\gamma)u_0 \\
  &=& \frac{n-2}{2}f u_0^\beta - \frac{\partial u_0}{\partial \eta}
 - \frac{n-2}{2}h_gu_0\\ 
 &\leq & 0
\end{eqnarray*}
since $u_0=u^+$ is a super-solution. By construction $u_1-u_0=0$ on $N$.
 
The maximum principle implies $u_1\leq u_0$ and the first step of the induction follows.

Assume, by induction, that $u_i\leq u_{i-1}$.
 
 Then, $L_\lambda (u_{i+1}-u_i)=\Delta u_{i+1}-\lambda u_{i+1} - \Delta u_i
 +\lambda u_i = -\lambda u_i +\lambda u_{i-1} = \lambda (u_{i-1}-u_i)
 \geq 0$.
 
 On $\partial M$ we have:
 
 \begin{eqnarray*}
 B_\gamma (u_{i+1}-u_i) & = & \frac{n-2}{2}fu_i^\beta + \gamma u_i
 -\frac{n-2}{2}fu_{i-1}^\beta -\gamma u_{i-1}\\
&= & \frac{n-2}{2} f(u_{i}^\beta-u_{i-1}^\beta)+\gamma (u_{i}-u_{i-1}).
\end{eqnarray*} 

If $f$ is nonnegative, then the above quantity is nonpositive
by induction hypothesis.

On the other hand, if there exists $x\in \partial M$ with $f(x)\leq 0$, then by choosing
$\gamma > \frac{n-2}{2}\beta  \|f\| \|u^+\|_{\partial M}$ we get 
\[\frac{n-2}{2}f(u_{i}^\beta-u_{i-1}^\beta)+\gamma (u_{i}-u_{i-1})
\leq 0,\] 
so the inequality $B_\gamma (u_{i+1}-u_i)\leq 0$ follows from 
the fact that 
\[u_{i-1}^\beta - u_i^\beta \leq \beta (u^+)^{\beta -1}(u_{i-1}-u_i).\]
Together with the fact that $u_{i+1}-u_i=0$ on $N$, it follows by the
maximum principle that $u_i$ is non-increasing.

We now show that $u_-\leq u_i$.

By hypothesis, $u_-\leq u^+=u_0$.
Assume, by induction, that $u_-\leq u_{i-1}$. Then $L_\lambda (u_--u_i)
=\Delta u_- - \lambda u_- -\Delta u_i+\lambda u_i = \Delta u_-
+\lambda (u_{i-1}-u_-) \geq 0$, by induction hypothesis and
the fact that $\Delta u_-\geq 0$. 

On $\partial M$ we have
\begin{eqnarray*}
B_\lambda (u_--u_i)& =& \frac{\partial u_-}{\partial \eta} + (\frac{n-2}{2}h_g
+\gamma)u_- - (\frac{n-2}{2}fu_{i-1}^\beta
+ \gamma u_{i-1})\\
&= & B (u_-)  
 + \frac{n-2}{2}f(u_-^\beta -u_{i-1}^\beta)+\gamma(u_--u_{i-1})\\
& \leq& \frac{n-2}{2}f(u_-^\beta -
u_{i-1}^\beta)+\gamma(u_--u_{i-1}).
\end{eqnarray*}
Should $f$ be 
positive, this last term would be non-positive by induction hypothesis.
Otherwise, $\gamma >\beta\frac{n-2}{2} \|u^+\|_{\partial B} \|f\|$ guarantees
$B_\lambda (u_--u_i)\leq 0$ since $u_k^\beta-u_-^\beta\leq \beta (u^+)^{\beta
-1}(u_k-u_-)$ for $k=1, \dots, i-1$. The fact that $u_i=u^+$ on $N$ and
$u_-\le u^+$ implies $u_--u_i\le 0$ on $N$. The claim
follows from the maximum principle.

The inequality $u_-\leq u_i\leq u^+$ in $M_1$ implies that the sequence
$u_i$ is uniformly bounded. From the first equation in (\ref{induct})
we conclude that $|\Delta u_i|$ is uniformly bounded as well.
Standard elliptic estimates imply that $\|u_i\|_{2,p}$ is uniformly
bounded for any $p>1$, and hence the Sobolev embedding implies
that there is a uniform bound for the sequence $u_i$ in the 
$C^{1,\nu}(\bar{M_1})$-norm. Differentiating the first equation in (\ref{induct})
we find that $|\nabla \Delta u_i|$ is uniformly bounded, and $L^p$ elliptic
estimates imply that $\|u_i\|_{3,p}$ is uniformly bounded for any $p>n$. The
compactness of the embedding $H^{3,p}(M_1)\hookrightarrow C^{2,\nu}(\bar{M_1}),
~0 <\nu <1-\frac{n}{p},~p>n,$ guarantees the existence of a subsequence
of functions $u_{i_k}$ converging to a function $u|_{M_1}\in C^{2,\nu}(\bar{M_1})$.
Because the sequence of functions $u_i$ is monotone we conclude that
the whole sequence converges to $u|_{M_1}$. That $u|_{M_1}$ is in $C^\infty(M_1)$
is a standard argument since it solves (\ref{cf}) on $M_1$.

 A diagonal procedure on an
exhaustion of $M$ by compact pieces like $M_1$ gives a way to construct
a globally defined smooth function $u\in C^\infty(M)$. Clearly $u_-\le u\le u^+$. 
Also, $u$ is a 
uniform limit of (a subsequence of) $u|{M_1}$'s over compact subsets, so 
it is straightforward to check that it is a solution 
to equation (\ref{cf}).
\end{proof}

%------------------------------------------------------------------------------------------------------------

\section{Existence of sub- and super-solutions} 
\label{E}

%------------------------------------------------------------------------------------------------------------

We construct an appropriate harmonic function that we will
use as a base for our sub- and super-solutions.

\begin{lemma} \label{4.1}
There exists $\mu >0$, and a positive smooth function 
$\mu \le v\le 1+\mu$ on $M$, with

\begin{eqnarray*}
\left\{
\begin{array}{clllr}
\Delta_g v &  = & 0 &\mbox{in}~M,\\
B_g v &  < & 0 & \mbox{on}~\partial M,\\
v &\sim  & 1+\mu & \mbox{near infinity}.
\end{array}
\right.
\end{eqnarray*}

\end{lemma}

\begin{proof}
Let $R>0$ be large. There always exists a positive solution $v_{R}$ 
of the homogeneous problem

\begin{eqnarray*}
\begin{array}{ll}
(P_R) &
\left\{
\begin{array}{clllr}
\Delta v_R &  = & 0 &\mbox{in}~\{x~:~d(x,\partial M)<R\},\\
 v_R &  = & 0 & \mbox{on}~\partial M,\\
v_R & = & 1 & \mbox{on}~\{x~:~d(x,\partial M)=R\}.
\end{array}
\right.
\end{array}
\end{eqnarray*}
A standard argument shows that as $R_i\to \infty$, the sequence 
$v_{R_i}$ converges uniformly on compact sets to a harmonic 
function $0\le v_\infty \le1$.

\it Claim.
\rm $v_\infty \sim 1$ on each end's infinity.

Let $E$ be an end and $0\le v_E\le 1$ be a harmonic barrier function 
that vanishes on the boundary of a large ball intersected with $E$ 
and is asymptotic to $1$ (See \cite{LT95}). By the maximum principle,  
$v_E$ is smaller or equal than $v_\infty$ . This way, $v_\infty$ is non-zero 
and asymptotic to $1$ on all ends.

We get that $B v_\infty = \partial /\partial \eta (v_\infty)$, but 
$\partial /\partial \eta (v_\infty)<0$ by Hopf's principle, since $v_\infty$
attains its minimum along the boundary (recall that $\eta$ is the 
outward-pointing normal of the boundary).\\
Pick $\mu>0$ so that $v:=v_\infty+\mu$ still satisfies $B v<0$. 
This way, $v \ge\mu>0$ and $v$ is asymptotic to $1+\mu$, as desired.
\end{proof}

\begin{prop}\label{fin}
For appropriately small constants $\epsilon,\delta>0$,
$u_-:=\epsilon v$ is a sub-solution, and $u^+:=\epsilon v+\delta$ is a super-solution.
\end{prop} 

\begin{proof}
Let $v$ be as before.
Note that, since $v$ is positive on the boundary, it makes sense to
write $\epsilon v= O(\epsilon)$ on $\partial M$. This way, for
 $\epsilon, \delta>0,~\beta>1$, one
has 
\[ (\epsilon v+\delta)^\beta = O(\epsilon^\beta)+O(\delta^\beta) ~\mbox{on}~\partial M.\] 
By definition, $u_- \le u^+$, and both are harmonic.
In order for them to be sub- and super-solutions,
we just have to check their behavior on the boundary.

\it Claim 1. 
\rm 
$Bu_--\frac{n-2}{2} f (u_-)^\beta \le 0$ .

Recall that by construction, $Bv<0$ on the boundary. Hence
\begin{eqnarray*}
Bu_- -\frac{n-2}{2} f (u_-)^\beta & = & \epsilon B v -\frac{n-2}{2} f (\epsilon v)^\beta\\
&\le & -\epsilon\min_{\partial M} |B v|+\frac{n-2}{2} \max_{\partial M}|f| (\epsilon v)^\beta\\
&= & -O(\epsilon)+O(\epsilon^\beta)\\
&\le & 0
\end{eqnarray*}
by taking $\epsilon>0$ small enough.

\it Claim 2.
\rm $B u^+- \frac{n-2}{2}f(u^+ )^\beta\ge 0$.

We see that
\begin{eqnarray*}
B u^+ -\frac{n-2}{2} f(u^+ )^\beta & = & \epsilon B v-\frac{n-2}{2} f (\epsilon v+\delta)^\beta\\
&\ge& -\epsilon\max_{\partial M}|Bv|+\delta (\frac{n-2}{2}h_g)\\
&&- \frac{n-2}{2}\max_{\partial M}|f|(\epsilon v+\delta)^\beta\\
&  = & -O(\epsilon) +O(\delta)-O(\epsilon^\beta)-O(\delta^\beta).
\end{eqnarray*}
The above line can be made nonnegative
by choosing $\epsilon$ smaller than $\delta$,  and $\delta$ small (notice
the plus sign next to $O(\delta)$).
 
This way, $0<\mu\epsilon \le u_- \le u^+\le \epsilon+\mu \epsilon+\delta$ are 
sub- and super-solutions respectively.  
\end{proof}

\begin{proof}[Proof of Theorem \ref{thm1}]
 
The existence of $u$ satisfying (\ref{cf}) is granted by applying the above
Proposition \ref{fin} and Theorem \ref{main}. For the critical case, i.e. $\beta=n/(n-2)$. 
The completeness of the metric $\tilde{g}=u^{4/(n-2)}g$ follows from the lower
bound  $u\ge u_- \ge\mu\epsilon>0$.
\end{proof}

%------------------------------------------------------------------------------------------------------------

\appendix
\section{Construction of positive manifolds}

%------------------------------------------------------------------------------------------------------------

We show how to construct a large class
of noncompact complete positive manifolds with boundary. Basically,
these examples come from removing ''small'' submanifolds from positive 
{\it compact} manifolds with boundary. 

\begin{remark} Positivity of compact manifolds is equivalent
to positivity of the first eigenvalue of problem (\ref{CN}), since
completeness is not an issue. A compact manifold with
boundary is positive if and only if  its {\it Yamabe constant} is positive 
(see \cite{E92}).
\end{remark}

Let $(N^n,\bar{g}),~n\ge 3$, be a positive compact manifold with boundary. 
 Consider a
collection of submanifolds $\Sigma=\cup_{i=1}^{k}\Sigma_i^{n_i}$,
where each $\Sigma_i$ is a submanifold in the interior of $N$
of dimension $0\le n_i\le \frac{n-2}{2}$; put $\Sigma_i=\{p_i\}$ whenever
$n_i=0$.

We will construct a metric $g=u^{4/(n-2)}\bar{g}$ on $M=N\setminus \Sigma$,
that is complete, scalar flat and has positive mean curvature on the
boundary. Also, $(M,g)$ will have large ends and will remain positive.

For $p\in int(N)$ let $G_p>0$ denote the Green's function
for the conformal Laplacian on $(N,\bar{g})$, which always exists
and satisfies $L_{\bar{g}}G_p=\delta_p$
and $B_{\bar{g}}G_p=0$. This way, for $c>0,~G_p+c$
satisfies
\[L_{\bar{g}} (G_p+c)=\delta_p,~B_{\bar{g}}(G_p+c)=c\frac{n-2}{2}h_{\bar{g}}>0\]
since $(N,\bar{g})$ is positive.

By a construction on the Appendix of Schoen and Yau's paper \cite{SY79}
which involves the Green's function, one can find, for each $\Sigma_i$
of positive dimension, positive functions $G_i$ that are singular on
$\Sigma_i$ and satisfy $L_{\bar{g}}G_i=0$ on $N\setminus \Sigma_i$.

A simple argument like that of Proposition \ref{fin} shows that
for appropriate coefficients $a_i>0, c>0$ the function 
\[u=\sum_{\{i|n_i>0\}} a_iG_i +\sum_{\{i|n_i=0\}}a_{i}G_{p_i}+c\] 
is singular on $\Sigma$ and satisfies $L_{\bar{g}}u=0$ and $B_{\bar{g}}u>0$.
Therefore $g=u^{4/(n-2)}\bar{g}$ remains positive.
 
The large codimension of the $\Sigma_i$ guarantees, via the standard 
argument in \cite{SY79}, that the singularities of $u$ are strong enough to make 
$g=u^{4/(n-2)}\bar{g}$ complete with large ends.

%\bibliographystyle{alpha}

%\bibliography{2733.bib}
%\end{document}

\end{document}